\documentclass[11pt,english]{article}
\usepackage[T1]{fontenc}
\usepackage[latin9]{inputenc}
\usepackage{geometry}
\usepackage{url}
\usepackage{amsthm}
\usepackage{amsmath}
\usepackage{amssymb}

\theoremstyle{plain}
\theoremstyle{plain}
\newtheorem{thm}{Theorem}

\usepackage{babel}

\begin{document}

\title{A Combinatorial Result on Block Matrices}

\date{S. M. Sadegh Tabatabaei Yazdi and Serap A. Savari%
\thanks{The authors are with the Department of Electrical and Computer Engineering,
Texas A\&M University, College Station, TX 77843 USA. Their e-mail
addresses are \protect\url{sadegh@neo.tamu.edu} and \protect\url{savari@ece.tamu.edu}.%
}}
\maketitle
\begin{abstract}
Given a matrix with partitions of its rows and columns and entries
from a field, we give the necessary and sufficient conditions that
it has a non--singular submatrix with certain number of rows from
each row partition and certain number of columns from each column
partition.
\end{abstract}

\section{Introduction}

Let $G$ be a matrix with entries from a field. Also let $S$ and
$T$ respectively denote the set of indices of rows and columns of
matrix $G$ and $G\left(\mathbf{s},\mathbf{t}\right)$ for sets $\mathbf{s}\subseteq S$
and $\mathbf{t}\subseteq T$ denote the submatrix of $G$ which is
the intersection of rows with indices in $\mathbf{s}$ and columns
with indices in $\mathbf{t}.$ In this paper we show that:
\begin{thm}
\label{thm:matrix-transversal} Suppose that the row set is partitioned
to $S=S_{1}\cup\cdots\cup S_{m}$ and the column set is partitioned
to $T=T_{1}\cup\cdots\cup T_{n}$. Given non--negative integers $s_{1},\cdots,s_{m},t_{1},\cdots,t_{n}$,
with $\sum_{i=1}^{m}s_{i}=\sum_{j=1}^{n}t_{j}=R,$ there are disjoint
subsets $\mathbf{s}_{1}\subseteq S_{1},\cdots,\mathbf{s}_{m}\subseteq S_{m}$
and $\mathbf{t}_{1}\subseteq T_{1},\cdots,\mathbf{t}_{n}\subseteq T_{n}$
with $|\mathbf{s}_{i}|=s_{i}$ and $|\mathbf{t}_{j}|=t_{j},$ such
that $G(\bigcup_{i=1}^{m}\mathbf{s}_{i},\bigcup_{j=1}^{n}\mathbf{t}_{j})$
is nonsingular, if and only if for every $I\subseteq\left\{ 1,\cdots,m\right\} ,K\subseteq\left\{ 1,\cdots,n\right\} $
we have \[
\mbox{rank}(G(\bigcup_{i\in U}S_{i},\bigcup_{k\in K}T_{k}))\geq\sum_{i\in I}s_{i}+\sum_{k\in K}t_{k}-R.\]
 
\end{thm}
This result is an analogues of the Hall's mathcing theorem in graph
theory \cite{graph-theory-diestel} and its extension to matroids,
called Rado's theorem \cite{mirskey,welsh}. There are several extensions
to Rado's theorem to more general structures such as greedoids \cite{transversaltheorie,greedoids-1,greedoids-2,mirskey,perfect}.
However we are unaware of any structure that admits an extension of
Rado's theorem and includes matrices as an element. Our main tool
to prove Theorem \ref{thm:matrix-transversal} will be a generalization
of matroids on the sets to the product of sets called \emph{bimatroids}
or \emph{linking systems} \cite{bimatroids-kung,linking-schrijver}.
Matrices are also a special case of bimatroids. We will combine a
result of Kung on a matroid structure of bimatroid \cite{bimatroids-kung}
and the Rado's theorem to prove Theorem \ref{thm:matrix-transversal}.

\section{Proof of Theorem \ref{thm:matrix-transversal}}

First we introduce matroids:

Given a set $E$ and a function $r:2^{E}\rightarrow\mathbb{N}$ we
say that the pair $(E,r)$ defines a matroid if and only if:
\begin{enumerate}
\item $r(A)\le|A|$ for all $A\subseteq E$. 
\item If $A,B\subseteq E$ with $A\subseteq B$, then $r(A)\le r(B)$. 
\item For any $A,B\subseteq E$, we have $r(A\cup B)+r(A\cap B)\leq r(A)+r(B).$ 
\end{enumerate}
Kung \cite{bimatroids-kung} and Schrijver \cite{linking-schrijver}
have indpendently introduced an extension of matroids to a structure
over products of sets, which is called bimatroids or linking systems.
Given sets $X$ and $Y$ and a function $\lambda:2^{X}\times2^{Y}\rightarrow\mathbb{N}$
we say that the triple $(X,Y,\lambda)$ defines a bimatroid if and
only if: 
\begin{enumerate}
\item $\lambda(U,V)\leq\min\left\{ |U|,|V|\right\} $ for all $U\subseteq X$
and $V\subseteq Y$.
\item If $U\subseteq X$ and $V\subseteq Y$ then for every $U'\subseteq U$
and $V'\subseteq V$, $\lambda(U',V')\leq\lambda(U,V).$
\item For any $U_{1},U_{2}\subseteq X$ and $V_{1},V_{2}\subseteq Y$, we
have $\lambda(U_{1}\cap U_{2},V_{1}\cup V_{2})+\lambda(U_{1}\cup U_{2},V_{1}\cap V_{2})\leq\lambda(U_{1},V_{1})+\lambda(U_{2},V_{2})$.
\end{enumerate}
One classic example of bimatroids is the matrix $G,$ with set of
indices of rows $S$ and set of indices of columns $T$. Then if $\lambda(\mathbf{s},\mathbf{t})=\mbox{rank}(G(\mathbf{s},\mathbf{t}))$
for any $\mathbf{s}\subseteq S$ and $\mathbf{t}\subseteq T,$ the
triple $(S,T,\lambda)$ is a bimatroid \cite{bimatroids-kung,linking-schrijver}. 

Next we consider the set $E=S\cup T$ which is the disjoint union
of the indices of rows and columns of the matrix $G$. We also consider
the following rank function defined on the subsets of $E.$ For every
$\mathbf{s}\subseteq S$ and every $\mathbf{t}\subseteq T$ we let\begin{equation}
r(\mathbf{s}\cup\mathbf{t})=\mbox{rank}(G(\mathbf{s},T\backslash\mathbf{t}))+|\mathbf{t}|.\label{eq:rank-definition}\end{equation}
Kung \cite{bimatroids-kung} shows that $(E,r)$ is a matroid. For
the sake of completeness we prove this result here. Namely, we prove
that function $r$ satisfies the properties 1--3 in the definition
of a matroid. First property follows because $\mbox{rank}(G(\mathbf{s},T\backslash\mathbf{t}))\leq|\mathbf{s}|$
and hence, $r(\mathbf{s}\cup\mathbf{t})=\mbox{rank}(G(\mathbf{s},T\backslash\mathbf{t}))+|\mathbf{t}|\leq|\mathbf{s}|+|\mathbf{t}|=|\mathbf{s}\cup\mathbf{t}|.$
For the second property suppose that $\mathbf{s}_{1}\subseteq\mathbf{s}_{2}\subseteq S$
and $\mathbf{t}_{1}\subseteq\mathbf{t}_{2}\subseteq T$ then we need
to show that $r(\mathbf{s}_{1}\cup\mathbf{t}_{1})\leq r(\mathbf{s}_{2}\cup\mathbf{t}_{2})$
or,\[
\mbox{rank}(G(\mathbf{s}_{1},T\backslash\mathbf{t}_{1}))+|\mathbf{t}_{1}|\leq\mbox{rank}(G(\mathbf{s}_{2},T\backslash\mathbf{t}_{2}))+|\mathbf{t}_{2}|.\]
We have \begin{align*}
\mbox{rank}(G(\mathbf{s}_{1},T\backslash\mathbf{t}_{1}))+|\mathbf{t}_{1}| & \overset{(a)}{\leq}\mbox{rank}(G(\mathbf{s}_{2},T\backslash\mathbf{t}_{1}))+|\mathbf{t}_{1}|\\
 & \overset{(b)}{\leq}\mbox{rank}(G(\mathbf{s}_{2},T\backslash\mathbf{t}_{2}))+|\mathbf{t}_{2}|\end{align*}
 where $(a)$ follows because the rank of a matrix is at least as
large as the rank of its submatrix and $(b)$ follows because matrix
$G(\mathbf{s}_{2},T\backslash\mathbf{t}_{2})$ is formed by removing
$|\mathbf{t}_{2}|-|\mathbf{t}_{1}|$ columns from $G(\mathbf{s}_{2},T\backslash\mathbf{t}_{1}).$
Since removing each column reduces the rank at most one unit, therefore
$\mbox{rank}(G(\mathbf{s}_{2},T\backslash\mathbf{t}_{1}))-\mbox{rank}(G(\mathbf{s}_{2},T\backslash\mathbf{t}_{2}))\leq|\mathbf{t}_{2}|-|\mathbf{t}_{1}|$.
Finally to prove that the third property holds for function $r,$
from the appendix for every $\mathbf{s}_{1},\mathbf{s}_{2}\subseteq S,\mathbf{t}_{1},\mathbf{t}_{2}\subseteq T$
we have \[
\mbox{rank}(G(\mathbf{s}_{1},T\backslash\mathbf{t}_{1}))+\mbox{rank}(G(\mathbf{s}_{2},T\backslash\mathbf{t}_{2}))\leq\mbox{rank}(G(\mathbf{s}_{1}\cap\mathbf{s}_{2},T\backslash(\mathbf{t}_{1}\cap\mathbf{t}_{2})))+\mbox{rank}(G(\mathbf{s}_{1}\cup\mathbf{s}_{2},T\backslash(\mathbf{t}_{1}\cup\mathbf{t}_{2}))).\]
 By adding the above relationship to $|\mathbf{t}_{1}|+|\mathbf{t}_{2}|=|\mathbf{t}_{1}\cap\mathbf{t}_{2}|+|\mathbf{t}_{1}\cup\mathbf{t}_{2}|$
we find that property 3, holds for function $r$ and every $A=\mathbf{s}_{1}\cup\mathbf{t}_{1},B=\mathbf{s}_{2}\cup\mathbf{t}_{2}$
where $\mathbf{s}_{1},\mathbf{s}_{2}\subseteq S,\mathbf{t}_{1},\mathbf{t}_{2}\subseteq T.$
Therefore $(S\cup T,r)$ is a matroid.

We now use the Rado--Hall theorem for matroids to prove our theorem.
In a matroid $(E,r),$ a set $A\subseteq E$ is an independent set
if $r(A)=|A|.$
\begin{thm}
(Rado--Hall) Let $(E,r)$ be a matroid and $A_{1},\cdots A_{n}\subseteq E.$
Given non--negative integers $\ell_{1},\cdots,\ell_{n},$ there exists
disjoint subsets $\mathbf{a}_{1}\subseteq A_{1},\cdots,\mathbf{a}_{n}\subseteq A_{n}$
with $|\mathbf{a}_{i}|=\ell_{i}$ and $\mathbf{a}_{1}\cup\cdots\cup\mathbf{a}_{n}$
an independent set, if and only if, for every subset $I\subseteq\left\{ 1,\cdots,n\right\} $
the following holds\[
r(\bigcup_{i\in I}A_{i})\geq\sum_{i\in I}\ell_{i}.\]
 
\end{thm}
Now consider the matroid $(S\cup T,r)$ where $r$ is defined in (\ref{eq:rank-definition}).
We next consider a partition of rows $S=S_{1}\cup\cdots\cup S_{m}$
and a partition of columns $T=T_{1}\cup\cdots\cup T_{n}$ and the
partition of $S\cup T=S_{1}\cup\cdots\cup S_{m}\cup T_{1}\cup\cdots\cup T_{n}.$
Assign to each $S_{i}$ a non--negative integer $s_{i}$ and to each
$T_{i}$ a non--negative integer $t_{i}$. By Rado--Hall theorem,
there exist disjoint subsets $\mathbf{s}_{1}\subseteq S_{1},\cdots,\mathbf{s}_{m}\subseteq S_{m}$
and $\mathbf{t}'_{1}\subseteq T_{1},\cdots,\mathbf{t}'_{n}\subseteq T_{n}$
such that $|\mathbf{s}_{i}|=s_{i}$ and $|\mathbf{t}'_{i}|=|T_{i}|-t_{i}$
and $\bigcup_{i=1}^{m}\mathbf{s}_{i}\bigcup_{j=1}^{n}\mathbf{t}'_{j}$
is an independent set in $(S\cup T,r)$ if and only if for every $I\subseteq\left\{ 1,\cdots,m\right\} ,J\subseteq\left\{ 1,\cdots,n\right\} $
we have \begin{equation}
r(\bigcup_{i\in I}S_{i}\bigcup_{j\in J}T_{j})\geq\sum_{i\in I}s_{i}+\sum_{j\in J}(|T_{j}|-t_{j}).\label{eq:condition1}\end{equation}
Let $\mathbf{t}_{i}=T_{i}\backslash\mathbf{t}'_{i}$. We have: \begin{align*}
\mbox{rank}(G(\bigcup_{i=1}^{m}\mathbf{s}_{i},T\backslash(\bigcup_{j=1}^{n}\mathbf{t}'_{j})))+|\bigcup_{j=1}^{n}\mathbf{t}'_{j}| & \overset{(a)}{=}\mbox{rank}(G(\bigcup_{i=1}^{m}\mathbf{s}_{i},\bigcup_{j=1}^{n}\mathbf{t}_{j}))+\sum_{j=1}^{n}(|T_{j}|-t_{i})\\
 & \overset{(b)}{=}\sum_{i=1}^{m}s_{i}+\sum_{j=1}^{n}(|T_{j}|-t_{i})\end{align*}
 where $(a)$ follows by definition and $(b)$ follows by the independence
of $\bigcup_{i=1}^{m}\mathbf{s}_{i}\bigcup_{j=1}^{n}\mathbf{t}'_{j}$.
Therefore \begin{equation}
\mbox{rank}(G(\bigcup_{i=1}^{m}\mathbf{s}_{i},\bigcup_{j=1}^{n}\mathbf{t}_{j}))=\sum_{i=1}^{m}s_{i}.\label{eq:condition2}\end{equation}
 Also condition $\eqref{eq:condition1}$ can be rewritten as the following.
For every $I\subseteq\left\{ 1,\cdots,m\right\} ,J\subseteq\left\{ 1,\cdots,n\right\} $
we have,\begin{align*}
\mbox{rank}(G(\bigcup_{i\in U}S_{i},T\backslash(\bigcup_{j\in J}T_{j})))+|\bigcup_{j\in J}T_{j}| & \geq\sum_{i\in I}s_{i}+\sum_{j\in J}(|T_{j}|-t_{j}).\end{align*}
hence \[
\mbox{rank}(G(\bigcup_{i\in U}S_{i},T\backslash(\bigcup_{j\in J}T_{j})))\geq\sum_{i\in I}s_{i}-\sum_{j\in J}t_{j}.\]
If we let $K=\left\{ 1,\cdots,n\right\} \backslash J,$ then the above
condition is equivalent to\[
\mbox{rank}(G(\bigcup_{i\in U}S_{i},\bigcup_{k\in K}T_{k}))\geq\sum_{i\in I}s_{i}+\sum_{k\in K}t_{k}-\sum_{j=1}^{n}t_{j}.\]
 Let suppose that $\sum_{i=1}^{m}s_{i}=\sum_{j=1}^{n}t_{j}=R.$ Then
since $G(\bigcup_{i=1}^{m}\mathbf{s}_{i},\bigcup_{j=1}^{n}\mathbf{t}_{j})$
has $\sum_{i=1}^{m}s_{i}$ rows and $\sum_{i=1}^{m}s_{i}$ columns,
$\eqref{eq:condition2}$ is equivalent to asking for $G(\bigcup_{i=1}^{m}\mathbf{s}_{i},\bigcup_{j=1}^{n}\mathbf{t}_{j})$
to be nonsingular. Therefore there exists a non--singular matrix $G(\bigcup_{i=1}^{m}\mathbf{s}_{i},\bigcup_{j=1}^{n}\mathbf{t}_{j})$
if and only if, for every $I\subseteq\left\{ 1,\cdots,m\right\} ,K\subseteq\left\{ 1,\cdots,n\right\} $
we have \[
\mbox{rank}(G(\bigcup_{i\in U}S_{i},\bigcup_{k\in K}T_{k}))\geq\sum_{i\in I}s_{i}+\sum_{k\in K}t_{k}-R.\]

\end{document}